\documentclass{amsart}
\usepackage{amssymb,amsxtra}
\theoremstyle{plain}
\newtheorem{thm}{Theorem}[section]
\newtheorem{lem}[thm]{Lemma}
\newtheorem{cor}[thm]{Corollary}
\newtheorem{df}[thm]{Definition}
\newtheorem{pro}[thm]{Proposition}
\theoremstyle{remark}
\newtheorem{rem}[thm]{Remark}
\title[Cuntz-Pimsner algebras]{On certain Cuntz-Pimsner algebras}
\author[a.~kumjian]{Alex Kumjian}
\address{Department of Mathematics, University of Nevada, Reno NV
89557, USA}
\email{alex@unr.edu}
\thanks{Research partially supported by \textsc{nsf} grant DMS-9706982}
\date{28 August 2001}
\begin{document}
\newcommand{\al}{\alpha}
\newcommand{\bt}{\beta}
\newcommand{\gm}{\gamma}
\newcommand{\Gm}{\Gamma}
\newcommand{\ep}{\varepsilon}
\newcommand{\eps}{\epsilon}
\newcommand{\et}{\eta}
\newcommand{\dl}{\delta}
\newcommand{\Dl}{\Delta}
\newcommand{\io}{\iota}
\newcommand{\kp}{\kappa}
\newcommand{\lm}{\lambda}
\newcommand{\ph}{\varphi}
\newcommand{\sg}{\sigma}
\newcommand{\tht}{\theta}
\newcommand{\zt}{\zeta}
\newcommand{\bd}{\partial}
\newcommand{\st}{\stackrel}
\newcommand{\x}{\times}
\newcommand{\rt}{\rtimes}
\newcommand{\ot}{\otimes}
\newcommand{\op}{\oplus}
\newcommand{\bop}{\bigoplus}
\newcommand{\emb}{\hookrightarrow}
\newcommand{\ra}{\rightarrow}
\newcommand{\lra}{\longrightarrow}
\newcommand{\imp}{\Rightarrow}
\newcommand{\lb}{\langle}
\newcommand{\rb}{\rangle}
\newcommand{\sub}{\subset}
\newcommand{\nul}{\emptyset}
\newcommand{\iso}{\simeq}
\newcommand{\Ad}{{\rm Ad}\,}
\newcommand{\Aut}{{\rm Aut}\,}
\newcommand{\im}{{\rm Im}\,}
\newcommand{\Hf}{{\mathfrak H}}
\newcommand{\Kf}{{\mathfrak K}}
\newcommand{\Ec}{{\mathcal E}}
\newcommand{\Fc}{{\mathcal F}}
\newcommand{\Gc}{{\mathcal G}}
\newcommand{\Hc}{{\mathcal H}}
\newcommand{\Kc}{{\mathcal K}}
\newcommand{\Lc}{{\mathcal L}}
\newcommand{\Oc}{{\mathcal O}}
\newcommand{\Tc}{{\mathcal T}}
\newcommand{\C}{{\mathbb C}}
\newcommand{\N}{{\mathbb N}}
\newcommand{\Q}{{\mathbb Q}}
\newcommand{\R}{{\mathbb R}}
\newcommand{\T}{{\mathbb T}}
\newcommand{\Z}{{\mathbb Z}}
\newcommand{\wt}{\widetilde}
\newcommand{\wh}{\widehat}
\newcommand{\ol}{\overline}
\newcommand{\qd}{\hfill $\Box$\pp}
\newcommand{\pr}{{\it Proof}\,:\ }
\newcommand{\tw}{\textwidth}
\newcommand{\vs}{\vskip1ex}
\newcommand{\amp}{\begin{minipage}[t]{.18\tw}}
\newcommand{\bmp}{\begin{minipage}[t]{.9\tw}}
\newcommand{\emp}{\end{minipage}}
\newcommand{\fn}{\footnotesize}
\newcommand{\scr}{\scriptstyle}
\newcommand{\rbx}{\raisebox}
\newcommand{\nnd}{\noindent}
\begin{abstract} 
Let $A$ be a separable unital C*-algebra and let 
$\pi : A \ra \Lc(\Hf)$ be a faithful representation of $A$ on a
separable Hilbert space $\Hf$ such that 
$\pi(A) \cap \Kc(\Hf) = \{ 0 \}$. 
We show that $\Oc_E$, the Cuntz-Pimsner algebra associated to the Hilbert
$A$-bimodule $E = \Hf \ot_{\C} A$, is simple and purely infinite.  
If $A$ is nuclear and belongs to the bootstrap class to which the UCT
applies, then the same applies to $\Oc_E$.
Hence by the Kirchberg-Phillips Theorem the isomorphism class of
$\Oc_E$ only depends on the $K$-theory of $A$ and the class of the unit.
\end{abstract}
\maketitle


In his seminal paper \cite{pm}, Pimsner constructed a C*-algebra
$\Oc_E$ from a Hilbert bimodule over a C*-algebra $A$ as a quotient of
a concrete C*-algebra $\Tc_E$, an analogue of the Toeplitz algebra,
acting on the Fock space associated to $E$. 
There has recently been much interest in these Cuntz-Pimsner algebras
(or Cuntz-Krieger-Pimsner algebras), which generalize both crossed
products by $\Z$ and Cuntz-Krieger algebras, as well as the associated
Toeplitz algebras. 
The structure of these C*-algebras is not yet fully understood,
although considerable progress has been made.  For example,  
Pimsner found a six-term exact sequence for the $K$-theory of $\Oc_E$
which generalizes the Pimsner-Voiculescu exact sequence 
(see \cite[Theorem 4.8]{pm});
conditions for simplicity were found in \cite{sc2,ms,kpw1,dpz}
and for pure infiniteness in \cite{z}.


The purpose of the present note is to analyze the structure of
Cuntz-Pimsner algebras associated to a certain class of Hilbert
bimodules.  Let $A$ be a separable unital C*-algebra and let 
$\pi : A \ra \Lc(\Hf)$ be a faithful representation of $A$ on a
separable Hilbert space $\Hf$ such that 
$\pi(A) \cap \Kc(\Hf) = \{ 0 \}$.  Then $E = \Hf \ot_{\C} A$ is a
Hilbert bimodule over $A$ in a natural way. 
We show that $\Oc_E$ is separable, simple and purely infinite.
If $A$ is nuclear and in the bootstrap class, then the same holds for
$\Oc_E$ and thus by the Kirchberg-Phillips theorem the
isomorphism class of $\Oc_E$ is completely determined by the
$K$-theory of $A$ together with the class of the unit (since $\Oc_E$
is $KK$-equivalent to $A$). 

Many examples of Cuntz-Pimsner algebras found in the
literature arise from Hilbert bimodules which are finitely generated
and projective; in such cases the left action must consist entirely of
compact operators.  Our examples do not fall in this class; in fact,
the left action has trivial intersection with the compacts.  And this
has some interesting consequences:
$\Oc_E \cong \Tc_E$ (see \cite[Corollary 3.14]{pm}) and the natural
embedding $A \emb \Oc_E$ induces a $KK$-equivalence 
(see \cite[Corollary 4.5]{pm}). 

In \S \ref{pre} we review some basic facts concerning the
construction of $\Tc_E$ as operators on the Fock space of $E$
and the gauge action $\lm : \T \ra \Aut(\Tc_E)$.  
We assume that the left action of $A$ does not meet the compacts $\Kc(E)$
and identify $\Oc_E$ with $\Tc_E$. 
The fixed point algebra $\Fc_E$, the analogue of the AF-core of a
Cuntz-Krieger algebra, contains a canonical descending sequence of
essential ideals indexed by $\N$ with trivial intersection.  The crossed 
product $\Oc_E \rt_{\lm} \T$ has a similar collection of essential ideals
indexed by $\Z$ on which the dual group of automorphisms acts in a
natural way.  By Takesaki-Takai duality 
$$
\Oc_E \ot \Kc(L^2(\T)) \cong (\Oc_E \rt_{\lm} \T) \rt_{\wh\lm} \Z;
$$
hence, much of the structure of $\Oc_E$ is revealed through an
analysis of the double crossed product. 

In \S \ref{result} we show that if $E$ is the Hilbert bimodule over
$A$ associated to a representation as described above, then for every
nonzero positive element $d \in \Oc_E$ there is a $z \in \Oc_E$ so
that $z^*dz = 1$; it follows that $\Oc_E$ is simple and purely
infinite (see Theorem \ref{spi}).  The proof of this proceeds
through a sequence of lemmas and is patterned on the proof of
\cite[Theorem 2.1]{rr}, which is in turn based on a key lemma of Kishimoto 
(see \cite[Lemma 3.2]{ks}).  Our argument uses the version of this
lemma found in \cite[Lemma 7.1]{op3} and this requires that we show
that the Connes spectrum of the dual action is full (this is also an
ingredient in the proof of simplicity found in \cite{dpz}).  
We invoke a version of a key lemma of R\o rdam for crossed products by
$\Z$ which arise from automorphisms with full Connes spectrum. 
The fact that $\Oc_E$ embeds equivariantly into 
$(\Oc_E \rt_{\lm} \T) \rt_{\wh\lm} \Z$ allows us to apply this lemma
to $\Oc_E$.  In \S \ref{app} we use the Kirchberg-Phillips
theorem to collect some consequences of this theorem as indicated
above and discuss certain connections with reduced (amalgamated) free
products.  

We fix some notation and terminology.
Given a C*-algebra $B$ we let $\wh{B}$ denote its spectrum, that is,
the collection of irreducible representations modulo unitary
equivalence endowed with the Jacobson topology (see \cite[\S 4.1]{pd}).
If $I$ is an ideal in a C*-algebra $B$, then every irreducible
representation of $I$ extends uniquely to an irreducible
representation of $B$.  This allows one to identify $\wh{I}$ with an
open subset of $\wh{B}$, the complement of which consists of the 
classes of irreducible representations which vanish on $I$. 
Given a *-automorphism $\bt$ of a C*-algebra $B$, let $\Gm(\bt)$
denote the Connes spectrum of $\bt$ (see \cite{o,co} or 
\cite[\S 8.8]{pd}); recall that 
$$
\Gm(\bt) = \bigcap_H {\rm Sp}\,(\bt|_H)
$$
where the intersection is taken over all $\bt$-invariant hereditary
subalgebras $H$.  A C*-algebra is said to be purely infinite if every
hereditary subalgebra contains an infinite projection.

I wish to thank D.~Shlyakhtenko for certain helpful remarks relating
to material in \S \ref{app}.
\section{Preliminaries}\label{pre}

We review some basic facts concerning Cuntz-Pimsner algebras; we shall
be mainly interested in those which arise from bimodules for which the
left action has trivial intersection with the compacts (see Remark
\ref{oe=te}). 
Let $A$ be a C*-algebra.  

\begin{df} {\rm (see \cite{ri,ka,l})}
Let $E$ be a right $A$-module.  Then $E$ is said to be a (right)
pre-Hilbert $A$-module if it is equipped with an $A$-valued inner
product $\lb \cdot, \cdot \rb_A$   
which satisfies the following conditions for all  $\xi, \et, \zt \in E$,  
$s, t \in \C$, and $a \in A$:
\begin{itemize}
\item[i.]  $\lb \xi, s \et + t\zt \rb_A  = 
         s\lb \xi,  \et \rb_A  + t\lb \xi,  \zt \rb_A$
\item[ii.] $ \lb \xi,  \eta a\rb_A  =  \lb \xi,  \et \rb_A  a $
\item[iii.] $ \lb \et,  \xi \rb_A  =  \lb \xi,  \et \rb_A^*$
\item[iv.] $\lb \xi, \xi \rb_A  \ge 0$ and  $\lb \xi, \xi \rb_A  = 0$ 
only if $\xi = 0$.
\end{itemize}
$E$ is said to be a (right) Hilbert $A$-module if it is complete in
the norm: $\|\xi\| = \| \lb \xi, \xi \rb_A \|^{1/2}$. 
\end{df}

Let $E$ be a Hilbert $A$-module.  Then $E$ is said to be {\em full} if
the span of the values of the inner product is dense.   
The collection of bounded adjointable operators on $E$, $\Lc(E)$, is a
C*-algebra. The closure of the span of operators of the form
$\tht_{\xi,\et}$ for $\xi, \et \in E$ 
(where $\tht_{\xi,\et}(\zt) = \xi\lb \et, \zt\rb_A$ for $\zt \in E$) 
forms an essential ideal in $\Lc(E)$ which is denoted $\Kc(E)$. A
Hilbert space is a Hilbert module over $\C$. 

\begin{df}
Let $E$ be a Hilbert $A$-module and $\ph : A \ra \Lc(E)$ be an
injective $*$-homomorphism.  Then the pair $(E, \ph)$ is said to be
Hilbert bimodule over $A$ (or Hilbert $A$-bimodule). 
\end{df}

Pimsner defines the Cuntz-Pimsner algebra $\Oc_E$ as a quotient of the
analogue of the Toeplitz algebra, $\Tc_E$, generated by creation
operators on the Fock space of $E$ (see \cite{pm}).  The injectivity
of $\ph$ is not really necessary (see \cite[Remark 1.2.1]{pm}).  We will
henceforth assume that $E$ is full (see \cite[Remark 1.2.3]{pm}).  

The Fock space of $E$ is the Hilbert $A$-module
$$
\Ec_+ = \bigoplus_{n = 0}^\infty E^{\ot n}
$$
where $E^{\ot 0} = A$ and for $n > 0$, $E^{\ot n}$ is the $n$-fold
tensor product: 
$$
E^{\ot n} = E \ot_A \cdots \ot_A E.
$$
Note that $\Ec_+$ is also a Hilbert $A$-bimodule with left action defined
by $\ph_+(a)b = ab$ for $a, b \in A = E^{\ot 0}$ and
$$
\ph_+(a)(\xi_1 \ot \cdots \ot \xi_n) =
\ph(a)\xi_1 \ot \cdots \ot \xi_n
$$
for $a \in A$ and $\xi_1 \ot \cdots \ot \xi_n \in E^{\ot n}$.

Then $\Tc_E \sub \Lc(\Ec_+)$ is the C*-algebra generated by the creation
operators $T_\xi$ for $\xi \in E$ where $T_\xi(a) = \xi a$ and
$$
T_\xi(\xi_1 \ot \cdots \ot \xi_n) = 
\xi \ot \xi_1 \ot \cdots \ot \xi_n.
$$
Observe that ${T_\xi}^*T_\et = \ph_+(\lb \xi, \eta\rb_A)$
for $\xi, \et \in E$.
Since $E$ is full, $\ph_+(A) \sub \Tc_E$; let $\io : A \emb \Tc_E$
denote the embedding.  Note that one may define $T_\xi$ for 
$\xi \in E^{\ot n}$ in an analogous manner and that we have  
${T_\xi}^*T_\et = \io(\lb \xi, \et\rb_A)$ 
for $\xi, \et \in E^{\ot n}$.

There is an embedding $\io_n : \Kc(E^{\ot n}) \emb \Tc_E$ (identify
$\Kc(E^{\ot 0})$ with $A$), given for $n > 0$ by 
$\io_n(\tht_{\xi,\et}) = T_{\xi}{T_{\et}}^*$ 
for $\xi, \et \in E^{\ot n}$.
Note that such operators preserve the grading of $\Ec_+$ and that there is
an embedding $\Kc(E^{\ot n}) \emb \Lc(E^{\ot m})$ for $m \ge n$. 
Let $C_n$ denote the C*-subalgebra of $\Tc_E$ generated by operators
of the form $T_{\xi}{T_{\et}}^*$ for $\xi, \et \in E^{\ot k}$ with 
$k \le n$ (by convention $C_0 = \io(A)$).  Then the $C_n$ form an
ascending family of C*-subalgebras. 

\begin{rem}\label{oe=te}
Suppose $\ph(A) \cap \Kc(E) = \{ 0 \}$; then the natural map 
$C_n \ra \Lc(E^{\ot m})$ is an embedding for $m \ge n$.
By \cite[Corollary 3.14]{pm} $\Tc_E \cong \Oc_E$ and the inclusion $A
\emb \Oc_E$ induces a $KK$-equivalence (see \cite[Corollary
4.5]{pm}).  Under the isomorphism of $\Tc_E$ with $\Oc_E$,
$\ol{\cup_n C_n}$ is mapped to $\Fc_E$, the analog of the AF core of a
Cuntz-Krieger algebra.
\end{rem}

For the remainder of this section we shall tacitly assume that 
$\ph(A) \cap \Kc(E) = \{ 0 \}$ and identify $\Tc_E$ with $\Oc_E$.

\begin{pro}\label{jn}
For each $n \in \N$ the C*-subalgebra, $J_n$, generated by
$\io_k(\Kc(E^{\ot k}))$ for $k \ge n$ is an essential ideal in
$\Fc_E$. We obtain a descending sequence of ideals   
$$
J_0 \supset J_1 \supset J_2 \supset \cdots
$$
with $J_0 = \Fc_E$ and  
$\cap_n J_n = \{ 0 \}$.  Furthermore, 
$J_n/J_{n+1} \cong \Kc(E^{\ot n})$ (thus $J_n/J_{n+1}$ is strong
Morita equivalent to $A$) and the restriction of the quotient map
yields an isomorphism $C_n \cong \Fc_E/J_{n + 1}$. 
\end{pro}
\begin{proof}
Given $n \in \N$ it is clear that $J_n$ is an ideal (see
\cite[Definition 2.1]{pm}).  To see that $J_n$ is essential it suffices to
show that for every $m$ and nonzero element $c \in C_m$ there is an element 
$d \in \Kc(E^{\ot k})$ for some $k \ge n$ such that $c\io_k(d) \ne 0$.
Let $k$ be an integer with $k \ge \max(m, n)$; since the map from $C_m$ to 
$\Lc(E^{\ot k})$ is an embedding for $k \ge m$, $c\xi \ne 0$ for
some $\xi \in E^{\ot k}$.  Then $cT_{\xi}{T_{\xi}}^* \ne 0$ and we take 
$d = \tht_{\xi,\xi}$.

The $J_n$ form a descending sequence of ideals by construction.
Since $\ph(A) \cap \Kc(E) = \{ 0 \}$, $C_m \cap J_n = \{ 0 \}$ 
for $m < n$.  Hence, $\cap_n J_n = \{ 0 \}$, for $\Fc_E$ is the
inductive limit of the $C_m$.  Further, for each $n$ we have 
$$
J_n = \io_n(\Kc(E^{\ot n})) + J_{n + 1}\quad\text{and}\quad
\io_n(\Kc(E^{\ot n})) \cap J_{n + 1} = \{ 0 \};
$$ 
it follows that $J_n/J_{n+1} \cong \Kc(E^{\ot n})$.  Finally, since
$$
\Fc_E = C_n + J_{n + 1}\quad\text{and}\quad
C_n \cap J_{n + 1} = \{ 0 \},
$$
we have $C_n \cong \Fc_E/J_{n + 1}$. 
\end{proof}

There is a strongly continuous action 
$$
\lm : \T \ra \Aut(\Oc_E)
$$
such that $\lm_t(T_\xi) = tT_\xi$.  The fixed point algebra under this
action is $\Fc_E$ and we have a faithful conditional expectation
$P_E: \Oc_E \ra \Fc_E$ given by
$$
P_E(x) = \int_\T \lm_t(x)\,dt.
$$
Consider the spectral subspaces of $\Oc_E$ under this action: 
for $n \in \Z$
$$
(\Oc_E)_n = 
\{ x \in \Oc_E : \lm_t(x) = t^nx \text{ for all }t \in \T \}.
$$

\begin{rem}\label{spec}
Note that $(\Oc_E)_n$ is the closure of the span of
elements of the form $T_{\xi}{T_{\et}}^*$ where $\xi \in E^{\ot k}$
and $\et \in E^{\ot l}$ with $n = k - l$. 
For $n \ge 0$ and $x \in (\Oc_E)_n$ we have $x^*x \in \Fc_E$ and 
$xx^* \in J_n$.  We may regard $(\Oc_E)_n$ as an 
$J_n$-$\Fc_E$-equivalence bimodule (see \cite{ri}); hence, $J_n$ is
strong Morita equivalent to $\Fc_E$ for each $n \ge 0$. If we regard 
$(\Oc_E)_1$ as a Hilbert $\Fc_E$-bimodule we have
(cf.\ \cite[\S 2]{pm} and \cite[\S 1.4]{sc2})
$$    
 (\Oc_E)_1  \cong E \ot_A \Fc_E,
$$
where the isomorphism is implemented by the map 
$\xi \ot a \mapsto T_\xi a$.
The crossed product $\Oc_E \rt_{\lm} \T$ may be identified with the
closure of the subalgebra of $\Oc_E \ot \Kc(\ell^2(\Z))$ consisting of
finite sums of the form
$$
\sum x_{ij} \ot e_{ij} 
$$
where $e_{ij}$ are the standard rank one partial isometries in 
$\Kc(\ell^2(\Z))$ and $x_{ij} \in (\Oc_E)_{j - i}$.
\end{rem}

Let $\wh{\lm}: \Z \ra \Aut(\Oc_E \rt_{\lm} \T)$ denote the dual
automorphism group. 

\begin{pro}\label{id}
There is an embedding $\eps : \Fc_E \emb \Oc_E \rt_{\lm} \T$ onto a
corner and a collection of essential ideals $\{ I_n \}_{n \in \Z}$ 
in $\Oc_E \rt_{\lm} \T$ satisfying the following conditions:
\begin{itemize}
\item[i.] For all $n \in \Z$, $\Fc_E$ is strong Morita equivalent to
$I_n$ and $A$ is strong Morita equivalent to $I_n/I_{n+1}$.
\item[ii.]  For all $n \ge 0$, $\eps(J_n) = \eps(1)I_n\eps(1)$.
\item[iii.] $I_n \sub I_m$ if $m \le n$.
\item[iv.] $\cap_n I_n = \{ 0 \}$
\item[v.] $\ol{\cup_n I_n} = \Oc_E \rt_{\lm} \T$
\item[vi.] $\wh{\lm}_k(I_n) = I_{n+k}$
\end{itemize}
\end{pro}
\begin{proof}
We use the identification of $\Oc_E \rt_{\lm} \T$ with a C*-subalgebra of 
$\Oc_E \ot \Kc(\ell^2(\Z))$ given in Remark \ref{spec}.  For each $n$
let $I_n$ be the ideal generated by $p_n = 1 \ot e_{nn}$.  Since 
$\Fc_E = (\Oc_E)_0$, it follows that $\Fc_E$ is isormorphic to the
corner determined by $p_n$ and thus is strong Morita 
equivalent to $I_n$.  The desired embedding  
$\eps : \Fc_E \emb \Oc_E \rt_{\lm} \T$ is given by
$\eps(a) = a \ot e_{00}$.  

Given an element of the form $a_{mn} = x_{mn} \ot e_{mn}$ in 
$\Oc_E \rt_{\lm} \T$ with $m \le n$, we have 
$$
{a_{mn}}^*a_{mn} = {x_{mn}}^*x_{mn} \ot e_{nn}
\quad\text{ and }\quad
a_{mn}{a_{mn}}^* = x_{mn}{x_{mn}}^* \ot e_{mm}
$$ 
with $x_{mn}{x_{mn}}^* \in J_{n - m}$; since $p_n$ may be expressed as a
finite sum of elements of the form ${a_{mn}}^*a_{mn}$, it follows that 
$I_n \sub I_m$ and that $p_mI_np_m = J_{n - m} \ot e_{mm}$.  
Thus $\eps(J_n) = \eps(1)I_n\eps(1)$ for all $n \ge 0$. 
Assertion (vi) follows from the fact that 
$\wh{\lm}_k(p_n) = 1 \ot p_{n + k}$.  
The remaining assertions follow from Proposition \ref{jn}.
\end{proof}

\section{$\Oc_E$ is simple and purely infinite}\label{result}
Let $A$ be a separable unital C*-algebra and let 
$\pi : A \ra \Lc(\Hf)$ be a faithful representation of $A$ on a
separable Hilbert space $\Hf$.

\begin{pro}\label{bimod}
With  $A$ and $\pi : A \ra \Lc(\Hf)$ as above,
$$
E = \Hf \ot_{\C} A
$$
is a full Hilbert bimodule over $A$ under the operations
$$
\lb \xi \ot a, \et \ot b \rb_A = \lb \xi, \et \rb a^*b,
\qquad
\ph(a)(\xi \ot b) = \pi(a)\xi \ot b
$$
for all $\xi, \et \in \Hf$ and $a, b \in A$.
Moreover, if  $\pi(A) \cap \Kc(\Hf) = \{ 0 \}$, then
$\ph(A) \cap \Kc(E) = \{ 0 \}$ and $\Oc_E \cong \Tc_E$.
\end{pro}
\begin{proof}
Note that $E = \Hf \ot_{\C} A$ is the tensor product of the Hilbert
$A$-$\C$-bimodule $\Hf$ and the Hilbert $\C$-$A$-bimodule $A$ 
as defined by Rieffel in \cite{ri} (see also \cite[Ch.\ 4]{l}).
The natural map from $\Lc(\Hf)$ to $\Lc(E) = \Lc(\Hf \ot_{\C} A)$ induces
an embedding $\Lc(\Hf)/\Kc(\Hf) \emb \Lc(E)/\Kc(E)$ (since $\Kc(\Hf)$
is mapped into $\Kc(E)$ and the Calkin algebra $\Lc(\Hf)/\Kc(\Hf)$ is
simple).  Hence, if $\pi(A) \cap \Kc(\Hf) = \{ 0 \}$, then  
$\ph(A) \cap \Kc(E) = \{ 0 \}$.  The last assertion, 
$\Oc_E \cong \Tc_E$, follows by \cite[Corollary 3.14]{pm}.
\end{proof}

Henceforth, we assume that $\pi(A) \cap \Kc(\Hf) = \{ 0 \}$ and
identify $\Oc_E$ with $\Tc_E$.
The aim of this section is to show that $\Oc_E$ is simple and purely
infinite.  Simplicity may be proven directly by invoking
\cite[Theorem 3.9]{sc2}: if $A$ is unital and $E$ is full, then $\Oc_E$ is
simple if and only if $E$ is minimal and nonperiodic.
Lemma \ref{hyp} would then be a consequence of \cite[Theorem 6.5]{op1}.
We follow a more indirect route patterned on the proof of 
\cite[Theorem 2.1]{rr}; this will also show that $\Oc_E$ is
purely infinite. 

\begin{rem}\label{rep}
With $E = \Hf \ot_{\C} A$ as above, we have 
$E^{\ot n} \cong \Hf^{\ot n} \ot_{\C} A$ via the map
$$
(\xi_1 \ot a_1) \ot (\xi_2 \ot a_2) \ot \cdots \ot (\xi_n \ot a_n) \mapsto 
(\xi_1 \ot \pi(a_1)\xi_2 \ot \cdots \ot \pi(a_{n-1})\xi_n) \ot a_n;
$$
similarly, if $\sg: A \ra \Lc(\Kf)$ is a representation of $A$ on a
Hilbert space $\Kf$, then  
$$
E^{\ot n}\ot_A \Kf \cong E^{\ot n-1} \ot_A \Hf \ot \Kf.
$$
Recall that the action of $\Fc_E$ on Fock space preserves the natural
grading.  Let $\wt\sg_n$ denote the representation of $\Fc_E$ on 
$E^{\ot n}\ot_A \Kf$ given by left action on $E^{\ot n}$.  Then the
restriction of $\wt\sg_n$ to $C_{n-1}$ is faithful: indeed, this
follows from the facts that the natural map  
$$
\Lc(E^{\ot n-1}) \ra \Lc(E^{\ot n-1} \ot_A \Hf \ot \Kf)
\cong \Lc(E^{\ot n}\ot_A \Kf)
$$
is an embedding (since $\pi$ is faithful) and that 
$\wt\sg_n|_{\Kc(E^{\ot n-1})}$ factors through $\Lc(E^{\ot n-1})$.
Note that $\wt\sg_n$ is equivalent to the representation of $\Fc_E$
obtained  from $\sg$ as follows: use the strong Morita 
equivalence between $A$ and $J_n/J_{n+1}$ 
to obtain a representation of $J_n/J_{n+1}$ and extend this to a
representation of $\Fc_E$.  Since the restriction of
$\wt\sg_n$ to $C_{n-1}$ is faithful, $\ker\wt\sg_n \sub J_{n}$ (see
Proposition \ref{jn}).  It follows that the closure of a point in 
$\wh{J_n} - \wh{J_{n+1}}$ contains the complement of $\wh{J_n}$.  
A similar assertion holds for $\Oc_E \rt_{\lm} \T$: 
for any $n \in \Z$ the closure of a point in $\wh{I_n} - \wh{I_{n+1}}$
contains the complement of $\wh{I_n}$. 
\end{rem}

\begin{lem}\label{hyp}
With $A$ and $E$ as above, $\Gm(\wh\lm_1) = \T$ where $\wh\lm$ is
the dual action of $\Z$ on $\Oc_E \rt_{\lm} \T$. 
\end{lem}
\begin{proof}
By \cite[Theorem 4.6]{op2} it suffices to find a dense invariant
subset of $(\Oc_E \rt_{\lm} \T)\sphat$\, on which ${\wh\lm_1}^*$ acts freely. 
That is, we must find an irreducible representation $\sg$ of 
$\Oc_E \rt_{\lm} \T$ such that, 
$\{ [\sg \circ \wh\lm_n] : n \in \Z \}$, the orbit of the
unitary equivalence class of $\sg$ under $\wh\lm^*$, is dense in 
$(\Oc_E \rt_{\lm} \T)\sphat$\, and 
$[\sg \circ \wh\lm_m] \ne [\sg \circ \wh\lm_n]$ if 
$m \ne n$.  Let $\sg_0$ be an irreducible representation of $A$ and
use the strong Morita equivalence between $A$ and $I_0/I_1$ to obtain
an irreducible representation $\sg'$ of $I_0/I_1$.  Then $\sg$, the
extension of $\sg'$ to $\Oc_E \rt_{\lm} \T$, is also irreducible. 
The classes $[\sg \circ \wh\lm_n]$ are distinct, for if $m < n$, 
$\sg \circ \wh\lm_m$ vanishes on $I_n$.
Moreover, for each $n \in \Z$ the closure of $[\sg\circ\wh\lm_n]$ in
$(\Oc_E \rt_{\lm} \T)\sphat$\, includes the classes of all irreducible
representations which vanish on $I_n$ (since
$[\sg\circ\wh\lm_n] \in \wh{I_n} - \wh{I_{n+1}}$,
see Remark \ref{rep}).  
Hence, $\{ [\sg \circ \wh\lm_n] : n \in \Z \}$ is dense in 
$(\Oc_E \rt_{\lm} \T)\sphat$.  
\end{proof}

Using Takesaki-Takai duality we show below that a C*-algebra $D$ equipped
with an action $\al$ of $\T$ may be embedded equivariantly as a corner in
$(D \rt_{\al} \T) \rt_{\wh\al} \Z$.  This fact is related to
Rosenberg's observation that the fixed point algebra under a compact
group action embeds as a corner in the crossed product (see \cite{ro}).

\begin{pro}\label{duality}
Given a unital C*-algebra $D$ and a strongly continuous action 
$\al : \T \ra \Aut(D)$,
there is an isomorphism $\psi$ of $D$ onto a full corner of 
$(D \rt_{\al} \T) \rt_{\wh\al} \Z$ which is equivariant in the sense
that $\wh{\wh{\al}}_t\circ\psi = \psi\circ\al_t$ for all $t \in \T$.
Moreover, $\psi(1) \in D \rt_{\al} \T$.
\end{pro}
\begin{proof}
By Takesaki-Takai duality \cite[7.9.3]{pd} there is an isomorphism
$$
\gm: D \ot \Kc(L^2(\T)) \cong (D \rt_{\al} \T) \rt_{\wh\al} \Z,
$$ 
which is equivariant with respect to $\al \ot \Ad \rho$ and
$\wh{\wh{\al}}$ (where $\rho$ is the right regular representation of
$\T$ on $L^2(\T)$).  The desired embedding is obtained by finding an
$\Ad\rho$ invariant minimal projection $p$ in $\Kc(L^2(\T))$ 
(cf.\ \cite{ro}):  set $\psi(d) = \gm(d \ot p)$ for $d \in D$.
Since $\psi$ is equivariant, $\psi(1)$ is in the
fixed point algebra of $\wh{\wh{\al}}$; hence, 
$\psi(1) \in D \rt_{\al} \T$.
\end{proof}

The following lemma is adapted from \cite[Lemma 2.4]{rr}; the proof is
patterned on R\o rdam's but we substitute 
\cite[Lemma 7.1]{op3} for \cite[Lemma 3.2]{ks}. 

\begin{lem}\label{key}
Let $B$ be a C*-algebra and let $\bt$ be an automorphism of $B$ such
that $\Gm(\bt) = \T$ and let $P$ denote the canonical conditional
expectation from $B \rt_\bt \Z$ to $B$.   Then for every positive
element $y \in B \rt_\bt \Z$ and $\ep > 0$ there are  
positive elements $x, b \in B$ such that 
$$
\| b \| > \| P(y) \| - \ep, \quad \| x \| \le 1 \quad \text{and} \quad
\| xyx - b \| < \ep.
$$  
If $y$ is in the corner determined by a projection $p \in B$, then
$x, b$ may also be chosen to be in the corner.
\end{lem}
\begin{proof}
As in the proof of \cite[Lemma 2.4]{rr} we may assume (by perturbing
$y$ if necessary) that $y$ is of the form 
$$
y = 
y_{-n}u^{-n} + \cdots + y_{-1}u^{-1} + y_0 + y_1u + \cdots + y_nu^n
$$
for some $n$ where $y_j \in B$ and $u$ is the canonical unitary in 
$B \rt_\bt \Z$ implementing the automorphism $\bt$; note that 
$y_0 = P(y)$ is positive.  

By \cite[Theorem 10.4]{op3} $\bt^k$ is properly outer for all 
$k \ne 0$.  Hence, by \cite[Lemma 7.1]{op3} there is a positive
element $x$ with $\| x \| = 1$ such that
$$
\| xy_0x \| > \| y_0 \| - \ep, \quad \text{and} \quad
\| xy_ku^kx \| = \| xy_k\bt^k(x) \| < \ep/2n
$$
for $0 < | k | \le n$. Set $b = xy_0x$; then a straightforward
calculation yields $\| xyx - b \| < \ep$.  
We now verify the last assertion.
Suppose that $y$ is in the corner determined by a projection 
$p \in B$; we may again assume that $y$ is of the above form.
Since $P$ is a conditional expectation onto $B$, 
$y_0 = P(y)$ is also in the corner determined by $p$. 
In the proof of \cite[Lemma 7.1]{op3} the positive element $x$ 
is constructed in the hereditary subalgebra determined by $y_0$;
hence we may assume that $x$ and therefore also $b = xy_0x$ lies
in the same corner. 
\end{proof}

Recall that $C_n$ is the C*-subalgebra of $\Fc_E$ generated by operators
of the form $T_{\xi}{T_{\et}}^*$ for $\xi, \et \in E^{\ot k}$ with 
$k \le n$ and that they form an ascending family of C*-subalgebras
with dense union. The subspace $E^{\ot n}$ is left invariant by $C_n$
and one has an embedding $C_n \emb \Lc(E^{\ot n})$.

\begin{lem}\label{easy}
Given a positive element $c \in C_n$ and $\ep > 0$, there is 
$\xi \in E^{\ot n}$ with $\| \xi \| = 1$ such that 
${T_\xi}^*cT_{\xi} \in C_0$ and 
$\| {T_\xi}^*cT_{\xi} \| > \| c \| - \ep$.
\end{lem}
\begin{proof}
The first assertion follows from a straightforward calculation: given  
$c \in C_n$ and $\xi \in E^{\ot n}$, then $c \xi \in E^{\ot n}$ and
$$
{T_\xi}^*cT_{\xi} = {T_\xi}^*T_{c\xi} = \io(\lb \xi, c\xi \rb_A) \in C_0.
$$
The second assertion follows from the embedding 
$C_n \emb \Lc(E^{\ot n})$ and the fact
$$
\| d \| = 
\sup\{ \| \lb \xi, d\xi \rb_A \| : \xi \in E^{\ot n}, \| \xi \| = 1 \} 
$$
for $d \in \Lc(E^{\ot n})$ positive.
\end{proof}

\begin{lem}\label{coup}
Given a positive element $a \in A$ and $\ep > 0$ with 
$\| a \| > \ep$, there is $\et \in E$ with 
$\| \et \| \le (\| a \| - \ep)^{-1/2}$ such that
${T_\et}^*\io(a)T_{\et} = 1$.
\end{lem}
\begin{proof}
Let $f$ be a continuous nonzero real-valued function supported on the
interval $[ \| a \| - \ep, \| a \| ]$ and choose a vector 
$\zt \in \pi(f(a))\Hf$ such that $\lb \zt, \pi(a)\zt \rb = 1$; 
we have 
$$
(\| a \| - \ep)\| \zt \|^2 \le \|\lb \zt, \pi(a)\zt \rb\| = 1.
$$
Then $\et = \zt \ot 1 \in E$ satisfies the desired conditions.
\end{proof}
It will now follow that $\Oc_E$ is simple and purely infinite (cf.\
proof of \cite[Theorem 2.1]{rr}). 

\begin{thm}\label{spi}
For every nonzero positive element $d \in \Oc_E$ there is a
$z \in \Oc_E$ so that $z^*dz = 1$.  Hence, $\Oc_E$ is simple and
purely infinite. 
\end{thm}
\begin{proof}
Let  $d \in \Oc_E$ be a nonzero positive element and choose $\ep$
so that $0 < \ep < \| P(d) \|/4$.  By Proposition \ref{duality} there
is a $\T$-equivariant isomorphism $\psi$ from $\Oc_E$ onto a corner of
$(\Oc_E \rt_{\lm} \T) \rt_{\wh\lm} \Z$ determined by a
projection $p \in \Oc_E \rt_{\lm} \T$.  We now apply Lemma \ref{key}
to the element $y = \psi(d)$ and the automorphism $\beta = \wh{\lm}_1$
(note $\Gm(\wh\lm_1) = \T$ by Lemma \ref{hyp}).  We identify
$\Oc_E$ with the  corner determined by $p$; note that under this
identification $\Fc_E$ is identified with $p(\Oc_E \rt_{\lm} \T)p$.
There are then positive elements $x, b \in \Fc_E$ so that 
$$
\| b \| > \| P(d) \| - \ep, \quad 
\| x \| \le 1 \quad \text{and} 
\quad \| xdx - b \| < \ep.
$$ 
Since $\cup_n C_n$ is dense in $\Fc_E$ we may
assume that $b \in C_n$ for some $n$. Hence, by Lemma \ref{easy} there
is $\xi \in E^{\ot n}$ with $\| \xi \| = 1$ such that 
$$
{T_\xi}^*bT_{\xi} \in C_0\quad \text{and} \quad 
\| {T_\xi}^*bT_{\xi} \| > \| b \| - \ep.
$$
Let $a$ denote the unique element of $A$ such that 
$\io(a) = {T_\xi}^*bT_{\xi}$; then $\| a \| > \| P(d) \| - 2\ep$ and 
$$
\| {T_\xi}^*xdxT_{\xi} - \io(a) \| =
 \| {T_\xi}^*(xdx - b)T_{\xi} \| < \ep.
$$
By Lemma \ref{coup} there is $\et \in E$
such that ${T_\et}^*\io(a)T_{\et} = 1$ and  
$$
\| \et \| \le (\| a \| - \ep)^{-1/2} 
< (\| P(d) \| - 3\ep)^{-1/2} < \ep^{-1/2}. 
$$
It follows that 
$$
\| {T_\et}^*{T_\xi}^*xdxT_{\xi}{T_\et} - 1 \| =
\| {T_\et}^*({T_\xi}^*xdxT_{\xi} - \io(a))T_\et \| \le
\| {T_\xi}^*xdxT_{\xi} - \io(a) \|(\ep^{-1/2})^2 < 1.
$$
Therefore, $c = {T_\et}^*{T_\xi}^*xdxT_{\xi}{T_\et}$ is an
invertible positive element and we take $z = xT_{\xi}{T_\et}c^{-1/2}$.
\end{proof}

\section{Applications and concluding remarks}\label{app}

We collect some applications of the above theorem and consider certain
connections with the theory of reduced (amalgamated) free product
C*-algebras. First we consider criteria under which the
Kirchberg-Phillips Theorem applies  
(see \cite[Theorem C]{kr}, \cite[Corollary 4.2.2]{ph}). 

\begin{thm}\label{main}
Let $A$ be a separable nuclear unital C*-algebra which belongs to the
bootstrap class to which the {\sc uct} applies (see \cite{rs}); let 
$\pi : A \ra \Lc(\Hf)$ be a faithful representation of $A$ on a
separable Hilbert space $\Hf$ such that 
$\pi(A) \cap \Kc(\Hf) = \{ 0 \}$ 
and let $E$ denote the Hilbert $A$-bimodule $\Hf \ot_{\C} A$.
Then $\Oc_E$ is a unital Kirchberg algebra (simple,
purely infinite, separable and nuclear) which belongs to the bootstrap
class.  
Hence, the Kirchberg-Phillips Theorem applies and the isomorphism
class of $\Oc_E$ only depends on $(K_*(A), [1_A])$ and not on the
choice of representation $\pi$.  
\end{thm}
\begin{proof}
First note that $\Oc_E$ is simple and purely infinite by Theorem
\ref{spi}. If $A$ is nuclear, then the argument given in the proof of
\cite[Theorem 2.1]{ds} shows that $\Oc_E$ must also be nuclear
(alternatively, the nuclearity of $\Oc_E$ follows from the
structural results discussed in \S \ref{pre}).  Hence, $\Oc_E$ is a
unital Kirchberg algebra.  Recall that the inclusion $A \emb \Oc_E$
defines a $KK$-equivalence (see \cite[Corollary 4.5]{pm}) which induces
a unit-preserving isomorphism $K_*(A) \cong K_*(\Oc_E)$.
Hence, if $A$ is in the bootstrap class, then $\Oc_E$ is also. 
Therefore, the Kirchberg-Phillips Theorem applies and the isomorphism
class of $\Oc_E$ only depends on $(K_*(A), [1_A])$. 
\end{proof}

Let $X$ be a second countable compact space, let $\mu$ be a nonatomic
Borel measure with full support and let 
$$
\pi : C(X) \ra \Lc(L^2(X, \mu))
$$
be the representation given by multiplication of functions.  Then
$\pi$ is faithful and
$$
\pi(C(X)) \cap \Kc(L^2(X, \mu)) = \{ 0 \}.
$$
Hence, we may apply the
above theorem with $A = C(X)$ and $\Hf = L^2(X, \mu)$.

\begin{cor}\label{abel}
Let $X$ and $\mu$ be as above.  Then
$$
E = L^2(X, \mu) \ot_{\C} C(X)
$$ 
is a Hilbert bimodule over $C(X)$ and $\Oc_E$ is a unital Kirchberg algebra.
The embedding  $C(X) \emb \Oc_E$ induces a (unit preserving) 
$KK$-equivalence.  Hence, the isomorphism class of $\Oc_E$ only depends on
$(K_*(C(X)), [1_{C(X)}])$ (and not on $\mu$); moreover, if $X$ is
contractible, then $\Oc_E \cong \Oc_\infty$.
\end{cor}

The following proposition is Theorem 5.6 of \cite{l}
(see also \cite[Theorem 3]{ka}); Lance calls this the 
Kasparov-Stinespring-Gelfand-Naimark-Segal construction.

\begin{pro}
Let $B$ and $C$ be C*-algebras, let $F$ be a Hilbert $C$-module and
let $f : B \ra \Lc(F)$ be a completely positive map, then there
is a Hilbert $C$-module $E_f$, a $*$-homomorphism  
$\ph_f : B \ra \Lc(E_f)$ and an element $v_f \in \Lc(F, E_f)$ such
that $f(b) = {v_f}^*\ph_f(b)v_f$ and $\ph_f(B)v_fF$ is dense in
$E_f$. 
\end{pro}

I am grateful to D.~Shlyakhtenko for the following observation.  Let
$\Tc$ denote the ``usual'' Toeplitz algebra (i.e.\ $\Tc_E$ where $E$
is the $1$-dimensional Hilbert bimodule over $\C$) and let $g$ denote
the vacuum state on $\Tc$.

\begin{pro}
Let $A$ be a separable unital C*-algebra and let 
$\pi : A \ra \Lc(\Hf)$ be a faithful representation of $A$ on a
separable Hilbert space $\Hf$ such that $\pi$ has a cyclic
vector $\xi \in \Hf$.  Let $f$ denote the vector state 
$\lb \xi, \pi(\cdot)\xi \rb$ and let $\tilde{f}$ denote the
corresponding completely positive map from $A$ to $\Lc(A)$ 
(given by $\tilde{f}(a) = f(a)1$).
Then $E = E_{\tilde{f}} \cong \Hf \ot A$ and $\Tc_{E}$ may be
realized as a reduced free product (see \cite{a,v}): 
$$
(\Tc_{E}, h) \cong (A, f) * (\Tc, g)
$$
for some state $h$ on $\Tc_{E}$.
\end{pro}
\begin{proof}
This follows from \cite[Theorem 2.3, Corollary 2.5]{sh}.
\end{proof}

As a result of this observation part (at least) of Corollary \ref{abel} follows
from the existing literature on reduced free products.  The simplicity
follows from a theorem of Dykema \cite[Theorem 2]{dy}.  Criteria for
when reduced free products are purely infinite have been found by
Choda, Dykema and R{\o}rdam in a series of papers \cite{dr1,dr2,dc};
but none seem to apply generally to the case considered in the corollary.

A theorem of Speicher (see \cite{sp}) on reduced amalgamated free
products (see \cite[\S 5]{v}) and Toeplitz algebras associated to Hilbert
bimodules yields a curious stability property of the algebras we have
been considering.  The following is the version given in 
\cite[Theorem 2.4]{bds}.
\begin{pro}
Suppose that $E_1$ and $E_2$ are full Hilbert bimodules over the
C*-algebra $A$.  Then
$$
\Tc_{E_1 \op E_2} = \Tc_{E_1} *_A \Tc_{E_2}.
$$
\end{pro}
We obtain the following corollary.  
\begin{cor}
Let $A$ be a separable nuclear unital C*-algebra which belongs to the
bootstrap class to which the {\sc uct} applies (see \cite{rs}) and let 
$\pi : A \ra \Lc(\Hf)$ be a faithful representation of $A$ on a
separable Hilbert space $\Hf$ such that 
$\pi(A) \cap \Kc(\Hf) = \{ 0 \}$.  
Let $E$ be the Hilbert bimodule $\Hf \ot_\C A$.  Then
$$
\Oc_E \cong \Oc_E *_A \Oc_E.
$$
\end{cor}
\begin{proof}
Observe that $E \op E = (\Hf \op \Hf) \ot_{\C} A$.  Since 
$\pi \op \pi : A \ra \Lc(\Hf \op \Hf)$ is a faithful representation and
$
(\pi \op \pi)(A) \cap \Kc(\Hf \op \Hf) = \{ 0 \},
$
the result follows follows from Theorem \ref{main} and the above
proposition. 
\end{proof}

\end{document}